\documentclass[10pt,draft]{amsart}

\usepackage{amssymb,latexsym,amsfonts,verbatim,amscd,ifthen}

\usepackage[usenames,dvipsnames]{pstricks}
\usepackage{pstricks}
\usepackage{epsfig}
\usepackage{pst-grad} 
\usepackage{pst-plot} 

\newtheorem{thm1}{Theorem}[section]
\newtheorem{lem1}[thm1]{Lemma}

\newtheorem{def1}[thm1]{Definition}

\newtheorem{cor1}[thm1]{Corollary}

\newtheorem{prop1}[thm1]{Proposition}
\newtheorem{ex1}[thm1]{Example}

\DeclareMathOperator{\HH}{\tilde{H}}
\DeclareMathOperator{\supp}{supp}
\DeclareMathOperator{\msupp}{m-supp}
\DeclareMathOperator{\vsupp}{v-supp}

\DeclareMathOperator{\Tor}{Tor}

 \DeclareMathOperator{\im}{Im}

\DeclareMathOperator{\bmax}{\textbf{max}\ }
 \DeclareMathOperator{\sign}{sign}
 \def\mapright#1{\smash{\mathop{\longrightarrow}\limits^{#1}}}
\newcommand\ds{\displaystyle}
\renewcommand\b{\beta}

\renewcommand\th{\theta}
\newcommand\D{\Delta}
\newcommand\tD{\widetilde{\Delta}}

\renewcommand{\L}{\mathcal{L}}
\newcommand{\A}{\mathcal{A}}

\begin{document}

\title[Syzygies]
{  On the generalized Scarf complex \\of lattice ideals }
\author[H. Charalambous]{Hara Charalambous}
\address { Department of Mathematics, Aristotle University of Thessaloniki,
Thessaloniki \\54124, GREECE
 } \email{hara@math.auth.gr}
\author[A. Thoma]{ Apostolos Thoma}
\address { Department of Mathematics, University of Ioannina,
Ioannina 45110, GREECE } \email{athoma@cc.uoi.gr}


\begin{abstract}
Let $\Bbbk$ be a field, $ \mathcal{L}\subset \mathbb{Z}^n$ be a
lattice  such that $\L\cap \mathbb{N}^n=\{ {\bf 0}\}$,
 and $I_\L\subset \Bbbk[x_1,\ldots, x_n]$ the corresponding lattice ideal.
We present the generalized Scarf complex of $I_\L$ and show that
it is indispensable in the sense that it is contained in every
minimal free resolution of $R/I_\L$.
\end{abstract}
\maketitle

\date{}

\section{Introduction}

\par Let $\Bbbk$ be a field,
  let $ \mathcal{L}\subset \mathbb{Z}^n$
be a lattice and $I_\L$ the corresponding lattice ideal in
$\Bbbk[x_1,\ldots, x_n]$. In a series of recent papers, see
 \cite{ATY}, \cite{BCMP}, \cite{CP}, \cite{ChKT},  \cite{DS},
\cite{Hi-O}, \cite{H-O} and \cite{AT},  the problem of getting
different minimal generating systems for lattice ideals $I_\L$ was
investigated. One of the motivations was a related  question from
Algebraic Statistics that asked when a prime lattice ideal
possesses a unique minimal system of binomial generators. As a
result the notion of indispensable binomials was introduced by
Ohsugi and Hibi \cite{Hi-O} to describe the binomials that up to
constant multiples are   part of all minimal systems of generators
of $I_\L$. The corresponding notion for higher syzygies is the
motivating question for the present article. Three landmark papers
dealing with free resolutions of lattice ideals  by Bayer, Peeva
and Sturmfels provide crucial leads: see \cite{BaSt}, \cite{PS1},
and \cite{PS2}. The algebraic Scarf complex is defined in
\cite{PS1} and it is shown to be  a minimal free resolution when
$I_\L$ is generic. The Taylor complex defined in \cite{BaSt} is
always a resolution, alas hardly ever a minimal one. The
complexity of computing minimal resolutions becomes even more
apparent in \cite{PS2}, and the importance of the topological
structure of the fibers is strongly hinted. Apart from the above,
the problem of computing syzygies of lattice ideals has also been
addressed by many recent articles, see for example \cite{BCMP2},
\cite{CG},
 \cite{P}, \cite{PVT}.

The structure of this paper is as follows: in section
\ref{intro_section}  we introduce the basic notions and
terminology. In section \ref{simplesyzygiessection} we introduce
simple minimal resolutions of $I_\L$:  these are complexes such
that the syzygies determined by
$\mathbb{Z}^n/\mathcal{L}$-homogeneous bases have {\em minimal}
support. In section \ref{generalizedScarf_section} we generalize
the  Scarf simplicial complex and basic fibers of \cite{PS1}. In
Section \ref{generalized_algebraic_Scarf_section} we introduce the
generalized algebraic Scarf complex and show that it is  an
indispensable complex. The last section contains examples and open
problems.

\section{Notation}
\label{intro_section}

Let $ \mathcal{L}\subset \mathbb{Z}^n$ be a lattice  such that
$\L\cap \mathbb{N}^n=\{ {\bf 0}\}$. The polynomial ring
$R=\Bbbk[x_1,\ldots, x_n]$ is positively multigraded by the group
$\mathbb{Z}^n/\mathcal{L}$, see \cite{MS}. Let ${\bf
a}_i=e_i+\mathcal{L}$ where $\{e_i:\  1\leq i\leq n\}$ is the
canonical basis of $\mathbb{Z}^n $. By $\mathcal{A}$ we denote the
subsemigroup of the group $\mathbb{Z}^n/\mathcal{L}$ generated by
$\{{\bf a}_i:1\leq i\leq n\}$. Since $\L\cap \mathbb{N}^n=\{ {\bf
0}\}$, the semigroup $\mathcal{A}$ is {\em pointed} meaning that
 $\{ x:\ x\in \mathcal{A} \textrm{ and } -x\in A\}=\{\bf{0}\}$.
Equivalently if $\A$ is a pointed semigroup generated by ${\bf
a}_1,\dots ,{\bf a}_n$ then  $ \mathcal{L} \subset \mathbb{Z}^n$
is the lattice of relations of ${\bf a}_1,\dots ,{\bf a}_n$ and
$\L\cap \mathbb{N}^n=\{ {\bf 0}\}$. We set $\deg_{\A}(x_i)={\bf
a}_i$.  The $ \mathcal{A}$-{\em degree} of the monomial ${\bf
x}^{{\bf v}}=x_1^{v_1} \cdots x_n^{v_n}$  is
\[ \deg_{ A}({\bf x}^{{\bf v}}):=v_1{\bf
a}_1+\cdots+v_n{\bf a}_n \in \mathcal{A}.\] When we want to put
the emphasis on $ \mathcal{L}$ we  occasionally write
$\deg_{\L}({\bf x}^{{\bf v}}):=\deg_{\A}({\bf x}^{{\bf v}})$. The
{\em lattice ideal} $I_{ \mathcal{L}}$ (or $I_\A$), associated to
$ \mathcal{L}$ is the ideal generated by all the binomials ${\bf
x}^{{\bf u}_+}- {\bf x}^{{\bf u}_-}$ where $\mathbf{u}_+,
\mathbf{u}_-\in \mathbb{N}^n$ and
$\mathbf{u}=\mathbf{u}_+-\mathbf{u}_- \in \mathcal{L}$.
 Prime lattice ideals are called toric ideals \cite{St} and  are the
defining ideals of toric varieties.
 For binomials in $I_{ \mathcal{L}}$, we define $\deg_{
\A}({\bf x}^{{\bf u}_+}- {\bf x}^{{\bf u}_-}):=\deg_{\A}{\bf
x}^{{\bf u}_+}$. Lattice ideals are $\A$-homogeneous. For $\bf{b}
\in \A$ we let $R[-{\bf b}]$ be the $\A$-graded free $R$-module of
rank 1 whose generator has $\A$-degree ${\bf b}$. Let
\[{(\bf F_\L,  \mathbf{\phi}}):
\quad 0\mapright{} F_p\mapright{\phi_p}\cdots
\cdots\mapright{}F_1\mapright{\phi_1} F_0\mapright{}
R/I_\L\mapright{} 0,\] be a minimal  $\A$-graded free resolution
of $R/I_L$.   The $i$-Betti number of $R/I_\L$ of $\A$-degree
$\mathbf{b}$ is equal to the rank of the $R$-summand of  $F_i$ of
$\A$-degree $\mathbf{b}$:
\[\beta_{i,{\mathbf{b}}}(R/I_L)=\dim_\Bbbk\Tor_i(R/I_\L,k)_{\mathbf{b}}\] and is denoted by
$\beta_{i,\mathbf{b}}(R/I_\L)$.  This is an invariant of $I_\L$,
see \cite{MS}.  For $\bf{b} \in \mathcal{A}$, we let $C_{\bf b}$
equal the fiber
\[C_{\bf
b}:=\deg_{\A}^{-1}(\bf b)=\deg_\mathcal{L}^{-1}({\bf b}):=\{x^u:\;
\deg_{\A}(x^u)=\bf b\}\] and
\[\D_{\bf b}:=\{F\subset \{1,\dots ,n\}:\ \exists\ x^{\bf a}\in C_{\bf b},\
F\subset \supp x^{\bf a}\}.\] It is well known that
\[ \b_{i,{\mathbf{b}}}(R/I_\L)=\dim_\Bbbk\HH_{i}(\D_{\bf b})\] see \cite{AH,  BCMP2, BH,
Stanley}. The degrees ${\bf b}$ for which $\HH_i(\D_{\bf b})\neq
0$ are called  {\em  i-Betti} degrees.

The semigroup $\A$ is pointed, so we can partially order $\A$ with
the relation
\[{\bf c} \geq {\bf d} \Longleftrightarrow \ \textrm{there is} \
{\bf e} \in \A \ \textrm{such that} \ {\bf c}={\bf d}+{\bf e}.\]
The minimal elements of the set $\{ {\bf b}:\ \beta_{i, {\bf b}
}(R/I_\L)\neq 0\}$ with respect to $\geq$ are called {\em minimal
$i$-Betti} degrees.

\section{Simple syzygies}
\label{simplesyzygiessection}

In this section we present the theory of simple syzygies of
$\A$-homogeneous ideals for a positive grading $\mathcal{A}$. We
will apply these results for lattice ideals.

Let $I$ be an $\A$-homogeneous ideal and $({\bf F},\phi)$ be a
minimal $\A$-graded free resolution of $R/I$, where
\[F_i=\ds{\bigoplus_{1\le t\le s_i}} R\cdot E_{ti}\ .\] In
particular $E_{10}$ is
 the basis element of $F_0\cong R$ of $\A$-degree 0. The
elements of $\im \phi_{i+1}=\ker \phi_i$ are called $i$-syzygies.
In the sequel all  syzygies are $\A$-homogeneous. We note that the
 zero syzygies of $R/I$ are the elements of $I$. Let $h$ be an
$\A$-homogeneous element of $F_i$. We write $h$ as a combination
of the basis elements $E_{ti}$ with nonzero coefficients:
\begin{equation}\label{syzygy}
h=\sum_{1\le t\le s_i} (\sum_{ c_{{\bf a_t}\neq 0}} c_{\bf a_t}
{\bf x}^{\bf a_t})E_{ti}\ .
\end{equation} Since $h$ is $\A$-homogeneous,
$\deg_\A{\bf x}^{\bf a_t}+\deg_\A(E_{ti})=\deg_\A(h)$. When $h\in
\ker \phi_{i}$ we  define $S(h)$, the {\em syzygy support} of $h$
to be the set
\[ S(h)=\{ {\bf x}^{\bf a_t}E_{ti}:\  c_{\bf a_t}\neq 0 \}, \]
and  partially order the $\A$-homogeneous $i$-syzygies by $h'\le
h$ if and only if $S(h')\subset S(h)$. We note that $E_{ti}$ may
appear
 in $S(h)$ more than once with different monomial coefficients.

\begin{def1}
 We say that
a nonzero $\A$-homogeneous $ h\in \ker \phi_i$ is {\it simple} if
there is no nonzero $\A$-homogeneous $ h'\in \ker \phi_i$ such
that $h'<h$.
\end{def1}

We note that the syzygy support of $h\in \ker \phi_{i}$ and thus
also the simplicity of $h$ depends on the basis $\{ E_{ti}\}$ of
$F_i$. We also note that when $I_\L$ is a lattice ideal then $h\in
I_\L$ is simple (as a zero syzygy of $R/I_L$) if and only if $h$
is a binomial. Next we define $\msupp(h)$, the {\em monomial
support} of $h\in F_i$. Let $C_1, C_2$ two fibers, $T_1\subset
C_1$, $T_2\subset C_2$. We let $T_1\cdot T_2=\{ m_1 m_2:\ m_i\in
T_i\}$. Let $h$ be a sum as in (\ref{syzygy}). We recursively
define $\msupp(h)$ by setting $\msupp(E_{10})=\{1\}$ and
\[\msupp(h)=\
\bigcup _{{\bf x}^{\bf a_t}E_{ti}\in S(h)}\ \{{\bf x}^{\bf
a_t}\}\cdot \msupp(\phi_i(E_{ti})).\]   We note that if the
$\A$-degree of $h$ is $\bf b$  then $\msupp(h)\subset C_{\bf b}$.
We also note that $\msupp(E_{ti})=\msupp(\phi_i(E_{ti}))$.

If $T$ is a subset of monomials then we set
\[\gcd(T):=\gcd (m:\ m\in T)\ .\]

\begin{def1} For a vector $\bf b\in \A$ we define the {\it gcd-complex}
$\D_{\gcd}(\bf b)$ to be the simplicial complex with vertices the
elements of the fiber $\deg_\mathcal{L}^{-1}(\bf b)$ and faces all
subsets $T\subset \deg_\mathcal{L}^{-1}(\bf b) $ such that
$\gcd(T)\neq 1$.
\end{def1}

Let $\bf b\in \A$. We remark that the gcd complex $\D_{\gcd}({\bf
b})$ and the complex $\D_{\bf b}$ have the same homology, see
\cite{ChThp}.

\begin{lem1}\label{simple_connected} Let $I$ be an $\A$-homogeneous ideal, $(\bf F, \phi)$
a minimal $\A$-graded free resolution of $R/I$ and  $h$ a simple
$i$-syzygy of $\A$-degree $\bf b$, $i\geq 1$. Then $\msupp(h)$ is
a connected subset of $\D_{\gcd}(\bf b)$.
\end{lem1}

\begin{proof} Let $h$ be given as in equation (\ref{syzygy}).
Since $\bf F_\L$ is  minimal it follows that ${\bf x}^{\bf
a_t}\neq 1$,  for all ${\bf x}^{\bf a_t}$ such that ${\bf x}^{\bf
a_t}E_{ti}\in S(h)$, see \cite{Stanley}. Therefore for any fixed
${\bf x}^{\bf a_t}$, the set $\{{\bf x}^{\bf a_t}\} \cdot
\msupp(\phi_i(E_{ti}))$ is a face of $\D_{\gcd}({\bf b})$.
Moreover, since $i\geq 1$  and $\bf F_\L$ is minimal, $1\notin
\msupp(\phi_i(E_{ti}))$. Therefore for each $E_{ti}$,
 the monomials in the sets $\{{\bf x}^{\bf a_t}\} \cdot
\msupp(\phi_i(E_{ti}))$ where ${\bf x}^{\bf a_t}E_{ti}\in S(h)$
are in the same connected component of $\D_{\gcd}({\bf b})$.
Suppose now that $\msupp(h)$ is disconnected and has $l\ge 2$
components. This means that there are disjoint index sets $J_r$,
$r=1,\ldots, l$  such that for each $r$
\[\bigcup_{j\in J_r}\ \  \bigcup_{{{\bf x}^{\bf a_j}E_{ji}\in S(h)}} \
\{{\bf x}^{\bf a_j}\} \cdot \msupp(\phi_i(E_{ji}))\] is a
connected component of $\msupp(h)$. Note that  the different
components of $\msupp(h)$ in  $\D_{\gcd}({\bf b})$ have no
variable in common. It follows that for a fixed $r$
\[h'=\sum_{j\in J_r} (\sum_{ c_{{\bf a_j}\neq 0}}c_{\bf a_j} {\bf
x}^{\bf a_j})E_{ji}\] is a syzygy and $h'< h$, a contradiction.
\end{proof}

We also remark that if $\bf (F, \phi)$ is a minimal $\A$-graded
free resolution of $R/I$ and  $\{H_{ti}: t=1,\ldots, b_i\}$ is an
$\A$-homogeneous basis for $F_i$ then by induction it can be shown
that $\gcd(\msupp(H_{ti}))=1$. Next we will consider minimal free
$\A$-graded resolutions of $R/I$ with a special property:

\begin{thm1}\label{simple_kernel} There exists a minimal $\A$-homogeneous generating set of
$\ker \phi_i$ consisting of simple i-syzygies.
\end{thm1}

\begin{proof}  Since  $\bf (F_\L, \phi)$ is a minimal $\A$-graded
free resolution of $S/I_\L$ and $\A$ is pointed  $\ker\phi_i$ can
be generated by $\A$-homogeneous elements. It is enough to show
that any $\A$-homogeneous  $i$-syzygy can be written as a sum of
simple $\A$-homogeneous $i$-syzygies. Let $h\in \ker\phi_i$ not
simple. We will use induction on $|S(h)|$. By hypothesis there
exists a simple $h'$ such that $h'<h$. Since $h'\not =0$  there is
an $\bf a$ such that
\[ h'={c_{\bf a}}'{\bf x}^{\bf a} E_{\bf a}+\sum_{{\bf b}\neq \bf
a}{c_{\bf b}}'{\bf x}^{\bf b}E_{\bf b},\quad  h=c_{\bf a}{\bf
x}^{\bf a} E_{\bf a}+\sum_{\bf b\neq \bf a}c_{\bf b}{\bf x}^{\bf
b}E_{\bf b},\quad c_{\bf a} {c_{\bf a}}'\neq 0.\] Note that if
${c_{\bf b}}'\not =0$ then $c_{\bf b}\not =0$. It follows that
\[ h=(h-\frac{c_{\bf a}}{{c_{\bf a}}'} h')+ \frac{c_{\bf a}}{{c_{\bf a}}'} h'\]
while $|S(h-\frac{c_{\bf a}}{{c_{\bf a}}'} h')|<|S(h)|$. Induction
finishes the proof.
\end{proof}

\begin{def1} Let $I$ be an $\A$-homogeneous ideal,
$(\bf{ F,\mathbf{\phi}})$ be an $\A$-graded free reso\-lution of
$R/I$ and for each $F_i$ we  let $B_i= \{E_{ti}\}$ be an
$\A$-homogeneous basis for $F_i$. We say that $\bf{ F}$ is {\em
simple} (with respect to the bases $B_i$) if $\phi_i(E_{ti})$ is
simple for each $i$ and $t$.
\end{def1}

The corollary follows immediately from Theorem
\ref{simple_kernel}.

\begin{cor1} Let $I$ be an $\A$-homogeneous ideal. There exists a
minimal simple $\A$-graded free resolution of $R/I$.
\end{cor1}

We note that there might be more than one minimal simple free
resolution of $R/I$.

\begin{ex1}\label{Koszul_example}{\rm  Let  $I=\langle f_1,\ldots, f_s \rangle$ be
a complete intersection lattice  ideal where $f_i: i=1,\ldots, s $
is an $R$-sequence of binomials and let
 $(\bf K_\bullet,\bf \theta)$ be the Koszul complex on the $f_i$.
The $\A$-homogeneous standard basis of  $\bf K_\bullet$ consists
of elements $E_J$ where $J $ ranges over all subsets of $[s]$; if
$J=\{k_1, \ldots, k_t\}$, where $k_1< \ldots < k_t$ then
\[ E_J=e_{k_1}\wedge \cdots\wedge e_{k_t},  \text{ and } \theta(E_J)=\sum
(-1)^{i+1} f_{k_i} E_{J\setminus \{k_i\}}.\] $(\bf K_\bullet,\bf
\theta)$ is a simple resolution of $R/I$ (with respect to the
basis $\{ E_J\}$). Indeed let  $h\not= 0$ be a syzygy such that
$h< \theta(E_J)$. Clearly this can only happen if $|J|>1$. Suppose
that $m_i E_{J\setminus \{k_i\}}\notin S(h)$ while $m_j
E_{J\setminus \{k_j\}}\in S(h)$; where $f_i=m_i-n_i, f_j=m_j-n_j$.
 Note that since $R$ has no zero divisors and the
coefficient of $E_{J\setminus \{k_i,k_j\}}$
 in $\theta(h)$  has to be zero, it follows that
 $n_i E_{J\setminus \{k_i\}}\in S(h)$. Therefore the coefficient of $E_{J\setminus \{k_i,k_j\}}$ in
$\theta(h)$ is $( c_in_i f_{k_j} - (c_jm_j-d_jn_j)f_{k_i})$ for
some $c_i, c_j, d_j \in \Bbbk$ which cannot be zero since
$f_{k_i}, f_{k_j}$ is an $R$-sequence, a contradiction.}
\end{ex1}

We note that  a simple syzygy is determined by its syzygy support:

\begin{thm1}\label{s_support_simple} If $h$ and $h'$ are two simple $i$-syzygies and $S(h)=S(h')$,
then there exists a $c\neq 0$ in $k$ such that $h=ch'$.
\end{thm1}

\begin{proof} Since $h'\neq 0$ there exists  ${c_{\bf a}}'\neq 0$  such that
$h'_t={c_{\bf a}}'{\bf x}^{\bf a}E_{\bf a}+\sum_{{\bf b}\neq {\bf
a}} c_{\bf b}{\bf x}^{\bf b}E_{\bf b}$. Therefore $c_{\bf a}\neq
0$ and
\[ h''=h-\frac{c_{\bf a}}{{c_{\bf a}}'} h' \in \ker \phi_i,\]
while $\ S(h'')\varsubsetneq S(h)$. Therefore $h''=0$ and
$h=c_{{\bf a}}/{c_{\bf a}}'h'$.
\end{proof}

Next we give the definition of an indispensable complex. We note
the related notion of rigidity, see \cite{M}.

\begin{def1}\label{indispensable} Let $I$ be an $\A$-homogeneous ideal.
We say that an $\A$-graded complex $(\bf G, \th)$ is an
indispensable complex for $R/I$ if for any simple minimal
$\A$-graded free resolution $(\bf F, \phi)$ of $R/I$ (with respect
to bases $B_j$ of $F_j$), there is an inclusion map  $i: (\bf G,
\theta)\mapright{} (\bf F, \phi)$ so that the image of $\bf G$ is
a subcomplex of $\bf F$. In particular for each $j$ there is a
subset ${B_j}'$ of $B_j$ so that the set $i^{-1}({B_j}')$ is an
$\A$-homogeneous basis of $G_j$. If $H_{tj}\in i^{-1}({B_j}')$ we
call $\th_j(H_{tj})$ an indispensable $j-1$ syzygy of $R/I$.

We say that $\bf G$ is a strongly indispensable complex for $R/I$
if the above holds without the requirement for $(\bf F, \phi)$ to
be simple.
\end{def1}

Let $I_\L$ be a toric ideal. A polynomial of $I_L$ is simple, as a
zero syzygy, if and only if it is a binomial. In this case the
indispensable  zero-syzygies of $R/I_\L$ are also strongly
indispensable.   More precisely in \cite{ChKT} it was shown that
the binomial $f$ of degree $\bf b$ is contained in every minimal
set of binomial generators of $I_\L$ if and only if $\bf b$ is a
minimal 1-Betti degree and the fiber $C_{\bf b}$ consists of just
two monomials with no common divisor: the difference of these two
monomials  is $f$ up to a constant multiple. Such a binomial is
not only contained in every minimal set of {\it binomial}
generators of $I_\L$; it is necessarily contained in {\it every}
minimal system of generators of $I_\L$ and thus is strongly
indispensable, see  \cite{PS1}. In \cite{PS1}, Peeva and Sturmfels
introduced the algebraic Scarf complex for $R/I_\L$. They  showed
that this complex is contained in every minimal resolution of
$R/I_\L$. It follows that the algebraic Scarf complex is a
strongly indispensable complex. We will generalize this
construction and show that the generalized algebraic Scarf complex
is also indispensable.

\section{ The generalized Scarf complex}
\label{generalizedScarf_section}

Let    $\bmax(\emptyset)=(-\infty,\ldots,-\infty)$. For  $J\subset
\mathcal{L}$, $0<|J|<\infty$,  we let
\[ \bmax(J)=(\max\{ { \bf a}_1:\ {{\bf a}\in J}\},\cdots  \ ,
\max\{{ \bf a}_n:\ {{\bf a}\in J}\}) \in \mathbb{Z}^n\ .\] We have
$\bmax(J')\le \bmax(J)$ if and only if $\bmax(J')_i\le \bmax(J)_i$
while $\bmax(\emptyset)$ is the smallest element. For $J\subset
K\subset \L$
 we define $\vsupp_K(J)$, the {\em variable support}
of   $J$ in $K$: if $J=\emptyset$ we set $
\vsupp_K(\emptyset)=\emptyset$ and for all other $J$ we set
\[\vsupp_K(J):=\{i:\ \exists \ {\bf a}\in J, \text{ such that }
\bmax(K)_i-{\bf a}_i>0\}.\] We note that
\[\vsupp_K(J)=\bigcup_{{\bf a}\in J}\supp({\bf x}^{\bmax(K)-{\bf
a}}),\] where $\supp({\bf x}^{{\bf a}})=\{i:x_i|{\bf x}^{{\bf
a}}\}$.

In \cite{PS1} the complex $\D_\L$ was defined to be the collection
of all finite subsets $J$ of $\L$ with unique $\bmax(J)$. We
extend this complex to  ${ \tD_\L}$ as follows:

\begin{def1} Let ${ \tD_\L}$ be the collection of all
finite subsets $J$ of $\L$ that satisfy the following conditions:
\begin{enumerate}
\item{} if $J'\subsetneq J$ then $\bmax (J')<\bmax( J)$ \item{}
          if ${\bf a}\notin J$ and $|J|\le 2$  then
                    ${\bf a}\nleq \bmax (J)$

                \item{} if ${\bf a}\notin J$, $|J|> 2$ and ${\bf a}\le \bmax (J)$
                then
                \[ \supp ({\bf x}^{\bmax (J)-{\bf a}})\bigcap \vsupp_J(J)=
                \emptyset\ .\]

\end{enumerate}
\end{def1}

\noindent We note that $\bmax(J)$ is determined by at most $n$
elements of $J$. If $J\in  \tD_\L$ then   the first condition
implies that $|J|\le n$. On the other hand if $|J|=2$ then the
first two conditions imply that $\bmax(J)$ is unique. Finally  all
sets $J$ with unique $\bmax (J)$ are in $\tD_\L$ and thus
$\D_\L\subset \tD_\L$.

\begin{prop1} $ \tD_\L$ is a simplicial complex.
\end{prop1}

\begin{proof} We have $\{{\bf a}\}\in \tD_\L$ for every $\bf a \in
\L$: ${\bf b} \le  {\bf a}\implies {\bf a}-{\bf b}\ge 0\in \L$, a
contradiction. The case $|J|=2$ is trivial. We now examine the
case $|J|>2, \ J\in \tD_\L$.

Let ${\bf a}\in J$; we will show that $J_1=J\setminus \{{\bf a}\}
\in \tD_\L$. Let $J_2\subsetneq J_1 $ such that $\bmax(J_2
)=\bmax(J_1)$. It follows that $\bmax(J_2\cup \{ {\bf
a}\})=\bmax(J)$, a contradiction.

Suppose now that ${\bf c}\notin J_1$ and ${\bf c}\le \bmax(J_1)$.
Therefore ${\bf c}\neq {\bf a}$, ${\bf c}< \bmax(J)$ and $\supp
({\bf x}^{\bmax (J)-{\bf c}})\bigcap \vsupp_J(J)=\emptyset$. On
the other hand since $\bmax(J_1)<\bmax(J)$, it follows that for
some $i$, $\bmax(J_1)_i<\bmax(J)_i$. Therefore ${\bf
a}_i=\bmax(J)_i$ while for all ${\bf b}\in J_1$, ${\bf
b}_i<\bmax(J)_i$. In particular $i\in \vsupp_J(J)$. This implies
that $i\notin \supp ({\bf x}^{\bmax (J)-{\bf c}})$, and ${\bf
c}_i=\max(J)_i$. Thus ${\bf c}\nleq \bmax(J_1)$, a contradiction.
It follows that $\bmax(J_1)$ is unique and $J_1\in \tD_\L$.
\end{proof}

There is a natural action of the lattice $\L$ in $\tD_\L$, since
$J\in  \tD_\L$ if and only if $J+{\bf a}\in  \tD_\L$ for any ${\bf
a}\in \L$. We identify $ \tD_\L$ with its poset of nonempty faces,
and we form the quotient poset ${ \tD_\L}/\L$. This poset is
called the {\it generalized Scarf complex} of $\L$.

\begin{prop1} The generalized Scarf complex ${ \tD_\L}/\L$ is a finite
poset.
\end{prop1}

\begin{proof}  Let   ${ \tD_\L}^0$ be the link to zero:
${ \tD_\L}^0=\{J\subset \L\setminus \{0\}: J\cup \{ 0\} \in
\tD_\L\}$. As in \cite{PS1}, since $\L$ acts transitively on the
vertices of $ \tD_\L$ it is enough to show that
  ${ \tD_\L}^0$
has finitely many vertices. The vertices $\bf a$ of ${ \tD_\L}^0$
are such that $\{ {\bf a}, {\bf 0}\}\in  \tD_\L$, therefore $\max(
\{{\bf a}, {\bf 0}\})={\bf a}^+$ is unique and we are exactly in
the case of \cite[Proposition 2.2]{PS1}:  there are finitely many
primitive elements of $\L$ and $\bf a$ is one of them, see also
\cite{BS}.
\end{proof}

For all $J\subset \L$ we define
\[
C_J:=\{ {\bf x}^{\bmax(J) -\bf a}:\ {\bf a}\in J\}\ .
\]
We note that $|C_J|=|J|$. Moreover if ${\bf x}^{\bf u}\in C_J$
then $\deg_\L {\bf x}^{\bf u}=\bmax(J)+L$. It follows that $C_J$
is a subset of the fiber $\deg^{-1}_\L (\bmax(J)+\L)$. We also
note that $C_J$ is not necessarily an entire fiber and there are
fibers or part of fibers that cannot be expressed in the form
$C_J$. The following Lemma determines exactly the cases when this
can happen.

\begin{lem1}\label{gcd}  Let $G$ be a subset of a fiber $C_{\bf b}$. Then
$G=C_J$ for some $J\subset \L$ if and only if $\gcd(G)=1$.
\end{lem1}

\begin{proof} Suppose that $G=C_J$ for some $J\subset \L$. For
each $1\leq i\leq n$ there exists ${\bf a} \in J$ such that
$\bmax(J)_i={\bf a}_i$. It follows that $i\notin \supp ( {\bf
x}^{\bmax(J)-{\bf a}})$. Therefore $\gcd(G)=1$.

Suppose that $\gcd(G)=1$ and  ${\bf x}^{\bf e}\in G$. Let $J=\{
{\bf a}\in \mathbb{Z}^n:\  {\bf x}^{{\bf e}-{\bf a}}\in G\}$.
Since $\deg_\L ( {\bf x}^{\bf e})= \deg_\L ( {\bf x}^{{\bf e}-{\bf
a}})=\bf b$ it follows that $J$ is a subset of $\L$. We claim that
$\bmax(J)={\bf e}$. Indeed, since ${\bf x}^{{\bf e}-{\bf a}}\in G$
it is clear that
  ${\bf e}\geq {\bf a}$ and therefore
$\bmax(J)\geq {\bf e}$. Moreover $\gcd (G)=1$ implies that  for
each
 $1\leq i\leq
n$ there exists an ${\bf a}\in J$ such that ${\bf a}_i={\bf e}_i$.
Therefore $\bmax(J)={\bf e}$ and  $G=C_J$.
\end{proof}

We isolate a slight variation of a useful remark of \cite{PS1}.

\begin{lem1}\label{hom_comp} If $\bf b$ is an $i$-Betti degree
then $C_{\bf b}=C_J$ for some subset $J$ of $\L$.
\end{lem1}

\begin{proof} Since $\HH_i(\Delta_{\bf b})\neq 0$
it is enough to show  that $\gcd(C_{\bf b})=1$, see Lemma
\ref{gcd}. Indeed if     $\gcd(C_{\bf b})\neq 1$ then $\Delta_{\bf
b}$ would be a cone with apex any variable in the support of
$\gcd(C_{\bf b})$ and $\Delta_{\bf b}$ would have no homology, a
contradiction.
\end{proof}

We will also make use of the following lemma:

\begin{lem1}\label{injectivity}  $C_J=C_{J'}$ if and only if $J'=u+J$, for some $u\in \L$.
\end{lem1}

\begin{proof} One direction is direct. For the opposite
let $a\in J$ and assume that $C_J=C_{J'}$. Then  there exists  $
a'\in J'$ such that ${\bf x}^{\bmax(J)-{\bf a}}= {\bf
x}^{\bmax(J')-{\bf a}'}$ and thus ${\bmax(J)-{\bf
a}}={\bmax(J')-{\bf a}'}$. We will show that $J=J'-{\bf a}'+{\bf
a}$. Let ${\bf b} \in J'-{\bf a}'+{\bf a}$. Then ${\bf b}+{\bf
a}'-{\bf a}\in J'$. We note that $C_{J'-a'+a}=C_{J'}=C_J$.
 Therefore \[{\bf
x}^{\bmax(J')-{\bf b}-{\bf a}'+{\bf a}}= {\bf x}^{( \bmax(J')-{\bf
a}')-{\bf b}+{\bf a}} = {\bf x}^{\bmax(J)-{\bf b}}\in C_J\] and
${\bf b}\in J$.
\end{proof}

When a fiber $C_{\bf b}$ is not of the form $C_J$ some of its
subsets may be expressed as such.

\begin{def1}\label{c-basic}
Let $C_{\bf b}$ be a fiber and $G\subset C_{\bf b}$. We say that
$G$ is a basic component of   $C_{\bf b}$ if the following are
satisfied:
\begin{itemize}
\item{} $G=C_J$ for some $J\subset \L$ \item{} $\gcd (G\setminus
\{m\})\neq 1$ for   all $m\in G$ and \item{} $G$ is a connected
component of $\D_{\gcd} ({\bf b})$ if $|C_{\bf b}|>2$.
\end{itemize}
If $C_{\bf b}$ satisfies the above properties then we call $C_{\bf
b}$ a basic fiber.
\end{def1}

\noindent We note that not all fibers contain subsets that are
basic fiber components. As a matter of fact we will show, see
Theorem \ref{finiteness_thm}, that the set of basic fiber
components is finite. The definition of a basic fiber first
appeared in \cite{PS1}. It follows directly from definition
\ref{c-basic} that if $C_{\bf b}$ is a basic fiber and $|C_{\bf
b}|>2$ then $C_{\bf b}$ has only one connected component. Moreover
if $|C_{\bf b}|=2$ then $C_{\bf b}$ is a basic fiber if and only
if $C_{\bf b}$ is disconnected. We also note that   if $G$ is a
basic component of $C_{\bf b}$ and $|C_{\bf b}|=2$ then $G=C_{\bf
b}$ is a basic fiber. This is the only way $C$ can be a basic
fiber component when $|C|=2$ as the following Lemma shows:

\begin{lem1}\label{card_2} If $|J|=2$ and $C_J$ is a basic component of $C_{\bf b}$ then
$C_J$ is a basic fiber.
\end{lem1}

\begin{proof} We have $|C_J|=2$.
By Lemma \ref{gcd} $C_J$ is disconnected in $\D_{\gcd} ({\bf b})$.
The third condition of definition \ref{c-basic} implies that
$C_J=C_{\bf b}$.
\end{proof}

In  \cite[Theorem 3.2, Lemma 3.3]{PS1} the following was shown:

\begin{lem1}\label{num_elements_basic_fiber} If $C_{\bf b}$ is a basic fiber with $i+1$ elements then
dim$\HH_i (\D_{\bf b})=1$ and $\bf b$ is a minimal $i$-Betti
degree.
\end{lem1}

Let $T$ be a subset of  a fiber $C_{\bf b}$. We denote by $[T]$
the set of monomials in $T$ divided by $\gcd(T)$. Let $I'=I\cup
\{m\}\subset T$, $m\notin I$. It is clear that
\[
[T \setminus I']= [ [T\setminus I]\setminus \{ \frac{m}{\gcd
(T\setminus I)}\}\  ].
\]

\begin{lem1}\label{subsets_c-basic} If $C_J$ is a basic component
of $C_{\bf b}$ and $\emptyset\neq I\subset C_J$ then $[
C_J\setminus I]$ is a basic fiber.
\end{lem1}

\begin{proof} By the preceding comment it suffices to prove the statement when
$|I|=1$. Let ${\bf d}=\deg_\A(\gcd (  C_J\setminus I))$ and ${\bf
b}'={\bf b} -{\bf d}$. Since    $\gcd  (  C_J\setminus I)\neq 1$
it follows that ${\bf b}'\neq {\bf b}$. We will show that
$[C_J\setminus I]=C_{{\bf b}'}$. Indeed if $m\in  C_{{\bf b}'}$
then $m {\bf x}^{\bf d}\in C_{\bf b}$ and clearly $m {\bf x}^{\bf
d}$ is in same connected component of $C_{\bf b}$ as $C_J\setminus
I$, thus $ m {\bf x}^{\bf d}\in C_J$ and   $[C_J\setminus
I]=C_{{\bf b}'}$. It is immediate that $\gcd [ C_J\setminus I]=1$.
The remaining condition follows as in  \cite[Lemma 2.4]{PS1}.
\end{proof}

The next lemma helps in computing $\gcd (C_J\setminus \{m\})$.

\begin{lem1}\label{gcd_subset_fiber} Let $C_J$
be a subset of $C_{\bf b}$, ${\bf a}\in J$ and $m=  {\bf
x}^{\bmax(J)- {\bf a}}$. Then
\[\gcd  (C_J\setminus \{m\}) =
  {\bf x}^{\bmax(J)-\bmax(J\setminus\{{\bf a}\})}\ .\]
\end{lem1}

\begin{proof}  An arbitrary element of  $C_J\setminus \{m\}$
is of the form ${\bf x}^{\bmax(J)-{\bf b}}$ where ${\bf b}\in
J\setminus \{ {\bf a}\}$. We note that ${\bf b}
 \le \bmax(J\setminus \{ {\bf a}\})$.
 It follows that
 \[
 {\bf x}^{\bmax(J)-{\bf b}}=
 {\bf x}^{\bmax(J)-\bmax(J\setminus\{ {\bf a}\})}
 {\bf x}^{(\bmax(J\setminus\{{\bf a}\})-{\bf b})}\ .
 \]
 For each $1\leq i\leq n$ there exist a ${\bf b}\in J\setminus \{ {\bf
 a}\}$ such that ${\bf b}_i=(\bmax(J\setminus\{{\bf a}\})_i$.
 Therefore $\gcd  (C_J\setminus \{m\})=
  {\bf x}^{\bmax(J)-\bmax(J\setminus\{a\})}$.
\end{proof}

The set of basic fiber components  forms a poset by setting
$C_{J'}\le C_{J}$ if and only if there exists a monomial ${\bf
x}^{\bf r}$ such that ${\bf x}^{\bf r} C_{J'}\subset C_J$. We
state and prove the analogue of \cite[Theorem 2.5]{PS1}.

\begin{thm1}\label{finiteness_thm}The poset of basic fiber components
 is isomorphic
to the generalized Scarf complex $\tD_\L/\L$ and is finite.
\end{thm1}

\begin{proof} Let $F$ be an element in ${ \tD_\L}/\L$. Choose a
representative $J$ of $F$. Let ${\bf b}=\bmax(J)+\L$. For any
other representative $J'$ of $F$, $C_J=C_{J'}\subset C_{\bf b}$.
We will show that $C_J$ is a basic component of $ C_{\bf b}$.

Let $m\in C_J$. We will show that $\gcd (C_J\setminus \{m\})\neq
1$. Suppose that $m={\bf x}^{\bmax(J)- {\bf a}}$, for ${\bf a} \in
J$. Since $J\in \tD_\L$, we have that $\bmax(J\setminus \{\bf
a\})<\bmax(J)$ and ${\bf x}^{\bmax(J)-\bmax(J\setminus\{a\})}\neq
1$. The desired inequality now follows by Lemma
\ref{gcd_subset_fiber}.

Let $|J|=2$. Then  $J\in \tD_\L$ implies  that $\bmax(J)$ is
unique and  $C_J= C_{\bf b}$. Suppose now that $|C_J|=|J|>2$.
Since $\gcd (C_J\setminus \{m\})\neq  1$, it follows that
$C_J\setminus \{m\}$ is a face of   $\D_{\gcd}({\bf b})$  and thus
$C_J$ is connected.  To show that $C_J$ is a connected component
of $C_{\bf b}$ we consider an element ${\bf x}^{\bf u}\in C_{\bf
b}\setminus C_J$. Then ${\bf u}+\L= \bf b$ and therefore ${\bf
u}=\bmax(J)-{\bf a}$ where ${\bf a}\in \L$ and ${\bf a}\notin J$.
Since ${\bf a}\le \bmax(J)$ and
\[ \supp ({\bf x}^{\bmax (J)-{\bf a}})\bigcap \vsupp_J(J)=
\emptyset\] we have that $\gcd({\bf x}^{\bf u},m)=1$, $\forall \
m\in C_J$. It follows that $C_J$ is a connected component of
$\D_{\gcd}({\bf b})$. Set $C_F=C_J$. It is immediate that  $\psi:
F\mapsto C_F$ is order preserving.

We will show that $\psi$ is bijective. The injectivity follows
immediately from Lemma \ref{injectivity}. To show that $\psi$ is
surjective we let $C_J$ be a basic component of $C_{\bf b}$. We
need to show that $J\in  \tD_\L$. Let $J'\subsetneq J$. For  $I=\{
{\bf x}^{\bmax(J)-{\bf a}}:\ {\bf a}\notin J'\}$ we have that
 $\gcd(C_J\setminus I)=\gcd(\{{\bf x}^{\bmax(J)-{\bf
a}}: {\bf a}\in J'\})\neq 1$. Therefore there is an $i$ such that
$\bmax(J)_i>{\bf a}_i$ for all ${\bf a}\in J'$. It follows that
${\bf a}<\bmax(J)$, $\forall \ {\bf a}\in J'$ and
$\bmax(J')<\bmax(J)$. Suppose  now that ${\bf a}\in \L$ is such
that ${\bf a}\notin J$ and ${\bf a}\le \bmax(J)$. It follows that
$\bmax(J\cup\{{\bf a}\})=\bmax(J)$ and $C_J\subsetneq C_{J\cup
\{a\}}\subset C_{\bf b}$. If $|J|=2$ then Lemma \ref{card_2} gives
a contradiction, so in this case $\bmax(J)$ is unique and $J\in
\tD_\L$. Suppose now that $|J|>2$. Since
 $m={\bf
x}^{\bmax(J)-{\bf a}}\in C_{\bf b}\setminus C_J$ and $C_J$ is a
connected component of $\D_{\gcd}({\bf b})$, it follows that no
variable in the support of $m$ is in the support of any monomial
in $C_J$ and thus
\[ \supp ({\bf x}^{\bmax (J)-{\bf a}})\bigcap \vsupp_J(J)=  \emptyset\ \]
as required.
\end{proof}

\section{The generalized algebraic Scarf complex}
\label{generalized_algebraic_Scarf_section}

We generalize the notion of the algebraic Scarf complex introduced
in \cite{PS1}.

\begin{def1} The generalized algebraic Scarf complex is the complex of
free $R$-modules
\[ ({\bf G_\L},\th)\ :=\bigoplus_{C\in { \tD_\L}/L} R\cdot E_{C}\]
where $E_{C}$ denotes a basis vector in homological degree $|C|-1$
and the sum runs over all basic fiber components $C$, identified
as elements of ${\tD_\L}/\L$. We let
\[ \th(E_{C})=\sum_{m\in C}\sign(m,C) \gcd(C\setminus\{m\})
E_{[C\setminus \{ m\}]},\] where $\sign(m,C)$ is $(-1)^{l+1}$ if
$m$ is in the $l^{th}$ position in the lexicographic ordering of
$C$.
\end{def1}

Our first remark is that $({\bf G_\L},\th)$ is a subcomplex of the
Taylor resolution, see \cite[Proposition 3.10]{BaSt}. Indeed the
canonical basis of the Taylor complex in homological degree $i$
consists of vectors $E_{C}$ where $C$ is a subset of a fiber  such
that $\gcd(C)=1$ and $|C|=i+1$. The differential of $({\bf
G_\L},\th)$ is the restriction of the differential of the Taylor
complex on the elements of $({\bf G_\L},\th)$.  Moreover we note
that the algebraic Scarf complex of \cite{PS1} is a subcomplex of
$({\bf G_\L},\th)$. Indeed the  canonical basis of the Scarf
complex in homological degree $i$ consists of vectors $E_{C}$
where $C$ is a basic fiber and the differential coincides for
these elements. We also note that for $i\le 1$ the algebraic Scarf
complex is identical to  the generalized algebraic Scarf complex.

\begin{thm1}\label{indispensable_theorem} The complex ${\bf G_\L}$ is an indispensable complex
for $R/I_\L$.
\end{thm1}

\begin{proof} Let $({\bf F_\L}, \phi)$ be a simple $\A$-graded minimal
resolution of $R/I_\L$ with respect to an $\A$-homogeneous basis
$\{ H_{tj},\ t=1,\ldots, b_j \}$  of $F_j$. We will use induction
on the homological degree $i$, the case $i=0$ being trivial. For
$i=1$, let $C=\{m_1, m_2\}$ be a basic fiber. Then
$\th_1(E_{C})=m_2-m_1$ is an indispensable binomial of $I_\L$
meaning that up to a constant multiple, $\th_1(E_{C})$ is part of
any minimal system of generators of the $I_\L$, see \cite{ChKT},
\cite{PS1}. Note that binomials in a lattice ideal are always
simple. Thus for $i\le 1$, $\th_i(E_{C})$ are indispensable.

Suppose now that ${\bf G_\L}$ is indispensable for homological
degrees less than $i$ and thus $\phi_j|_{G_j}=\theta_j$ for $j<i$.
Let $C$ be  a basic component of $C_{\bf b}$ with cardinality
$i+1$. It is clear  that $\th_i(E_{C})\in \ker \phi_{i-1}$. We
will show that $\th_i(E_{C})$ is a simple syzygy. The proof is
essentially the same as in Example \ref{Koszul_example}. Indeed
suppose that $0\neq h\in \ker \phi_{i-1}$ and $h< \th_i(E_{C})$.
Since $S(h)\subset S(\th_i(E_{C}))$ it follows that
 $h\in G_{i-1}$ and
$\phi_{i-1}(h)=\th_{i-1}(h)$. Moreover since $S(h)\subsetneq
S(\th_i(E_{C}))$ and $S(\th_i(E_{C}))=\{\ \gcd(C\setminus\{m\})
E_{[C\setminus \{ m\}]}:\ m\in C\}$, it follows that for some
monomial $m_1\in C$,
\[
\gcd(C\setminus\{m_1\}) E_{[C\setminus \{ m_1\}]} \notin S(h).
\]
On the other hand since $h\neq 0$ we can find $m_2\in C$ such that
\[ \gcd(C\setminus\{m_2\})E_{[C\setminus \{ m_2\}]} \in
S(h).\] It follows that $\gcd(C\setminus
\{m_1,m_2\})E_{[C\setminus \{m_1, m_2\}]}\in S(\phi_{i-1}(h))$, a
contradiction since $\phi_{i-1}(h)=0$.

Next suppose that $\th_i(E_{C})$ is an  $R$-linear combination of
$i$-syzygies of strictly smaller $\A$-degree. It follows that $\bf
F_\L$ in homological degree $i$ has a basis generator $h$ of
$\A$-degree $\bf b_1$, where ${\bf b_1}< \bf b$ and
\[
\msupp(\theta_i (E_{C}))\cap \msupp( mh) \neq \emptyset, \quad
\text{ for some } m\in C_{{\bf b}-{\bf b_1}}. \] By Lemma
\ref{hom_comp} we have that $\gcd(C_{\bf b_1})=1$. Therefore  $m
C_{\bf b_1}\cap C\neq \emptyset$. Since $C$ is a connected
component it follows that $m C_{\bf b_1}\subset C$ and therefore
 $C_{{\bf b_1}}=[C\setminus I]$ for some $I\subset C$.
Lemma \ref{subsets_c-basic} implies that $C_{\bf b_1}$ is basic
while Lemma \ref{num_elements_basic_fiber} implies that $C_{\bf
b_1}$ has at least $i+1$ elements. Recall that by construction
$|C|=i+1$. Since $\gcd(C)=1$, $m C_{J_1} \subset C$ and $|m C_{\bf
b_1}|=i+1$ it follows that $|C|>i+1$, a contradiction.

We still need to show that  we can identify  $E_{C}$ with a
constant multiple of a basis element $ H_{ti}$ for some $t$.
Suppose that  $\theta_i(E_{C})= \sum c_t \phi_i(H_{ti})$ where
$\deg_{\A} (H_{ti}) ={\bf b}$ for at least one $t$. Then  for some
$t$ with $\deg_{\A} (H_{ti}) ={\bf b}$, we have that
\[\gcd( C\setminus \{ m\})E_{[C\setminus \{ m\}]}\in
S(\theta_i(E_{C}))\cap S(\phi_i(H_{ti}))\ ,\] for a monomial $m\in
C$. Therefore $\msupp(\theta_i(E_{C}))\cap
\msupp(\phi_i(H_{ti}))\neq \emptyset$. \\Since
$\msupp(\theta_i(E_{C}))$ $=C$,
$\msupp(H_{ti})=\msupp(\phi_i(H_{ti}))$ and by Lemma
\ref{simple_connected} $\msupp(\phi_i(H_{ti}))$  is connected, it
follows that $\msupp(H_{ti})\subset C$. Suppose that
\[\phi_i(H_{ti})=\sum_s p_s H_{s,{i-1}} \]
where $H_{s,{i-1}}$ are basis generators of $F_{i-1}$ of
$\A$-degree ${\bf b_s}$ and $p_s\in R$ is $\A$-homo-geneous of
$\A$-degree ${\bf b}-{\bf b_s}$. It is clear that ${\bf b_s}<\bf
b$ and that $C_{\bf b_s} < C$. Moreover $\HH_{i-1}(\D_{\bf
b_s})\neq 0$ and therefore $C_{\bf b_s}=[C\setminus I]$ for some
subset $I$ of $C$. We note that the set $I$ need not be unique. By
Lemma \ref{subsets_c-basic} it follows that $C_{\bf b_s}$ is a
basic fiber. By induction $E_{C_{\bf b_s}}$ is indispensable and
is the unique basis element of $F_{i-1}$ of $\A$-degree ${\bf
b_s}$, up to a constant multiple. By Lemma
\ref{num_elements_basic_fiber} $|C_{\bf b_s}|\ge i$.   Since
$C_{\bf b_s}=[C\setminus I]$, it follows that $|C_{\bf b_s}|\le
i$. Therefore $|C_{\bf b_s}|= i$ and  $C_{\bf b_s}=[C\setminus \{
m_l\}]$ for a monomial $m_l$. Let $c{\bf x^{\bf \gamma}}$ be a
monomial term of $p_s$ where $c\in \Bbbk-\{0\}$. Since
$\msupp(\phi_i(H_{ti}))\subset C$, $\gcd(C)=1$ and $\gcd(C_{\bf
b_s})=1$ it follows that ${\bf x^{\bf \gamma}}=\gcd(C\setminus \{
m_l\})$, so actually $p_s=c{\bf x^{\bf \gamma}}$ for a $c\not= 0$.
Thus $\phi_i(H_{ti})\le \theta_i(E_{C})$ and since both are simple
we get that $S(\phi_i(H_{ti}))=S(\theta_i(E_{C}))$. We apply
Theorem \ref{s_support_simple} to obtain the desired conclusion.
\end{proof}

 Next we consider a complex that sits between the
algebraic Scarf complex and the generalized algebraic Scarf
complex.

\begin{def1}\label{strongly_algebraic} Let
$B_i=\{ E_C:\ C \textrm{ basic component of } C_{\bf b},$
$|C|=i+1,$ ${\bf b}$  minimal  $i\textrm{-Betti degree},$
$\dim_\Bbbk\HH_i(\D_{\bf b})=1\}$. The generalized  strongly
algebraic Scarf complex ${\bf S}_\L$ is the subcomplex of $\bf
G_\L$ with basis  in homological degree $i$ the set $B_i$.
\end{def1}

We note that if $C$ is a basic fiber then $C$ satisfies the
conditions of Definition \ref{strongly_algebraic} and thus ${\bf
S}_\L$ contains the algebraic Scarf complex of \cite{PS1}. This
containment can be strict as Example \ref{strongly_algebraic_not
Scarf} shows. We point out that if $C$ is a basic component of
$C_{\bf b}$ then for $i>0$ the reduced $i$-homology group of the
corresponding component of  $\D_{\bf b}$ has dimension $1$.

\begin{thm1}\label{strongly_indispensable_theorem}
The complex ${\bf S}_{\L}$ is a strongly indispensable complex for
$R/I_\L$.
\end{thm1}

\begin{proof} Let
$C_{\bf b}$ be a fiber with a component satisfying the conditions
of Definition \ref{strongly_algebraic}. These conditions imply the
following for any $\A$-graded free resolution of $R/I_\L$, $({\bf
F_\L}, \phi)$, and  $\{ H_{ti}:\ t=1,\ldots, b_i\}$  an
$\A$-homogeneous basis of $F_i$: there is a unique basis element
$H_{ti}$ of $\A$-degree ${\bf b}$ and there is no generator
$H_{li}$ of smaller $\A$-degree. Following the proof of Theorem
\ref{indispensable_theorem} this automatically implies that
$\theta_i(E_{C})=\phi_i(H_{ti})$.
\end{proof}

In \cite{PS1} it was shown that whenever the ideal $I_\L$ is
generic, meaning that the support of each minimal binomial
generator of $I_\L$ is $[n]$, then the algebraic Scarf complex is
a minimal resolution of $R/I_\L$. The converse is not necessarily
true as example \ref{Scarf_resolution_not_generic} shows: there
are ideals which are not generic but the algebraic Scarf complex
is a minimal resolution of $R/I_\L$. We note the following:

\begin{thm1}If $I_{\L}$ is generated by indispensable binomials
then the generalized  Scarf complex equals the  Scarf complex.
\end{thm1}

\begin{proof} Suppose that $C_{\bf b}$ is not a basic fiber so
that $C\neq C_{\bf b}$ is a basic fiber component. It is immediate
that  $C_{\bf b}$ has more than one connected components. By
\cite{ChKT} any binomial  which is  the difference of two
monomials belonging to different connected components of $C_{\bf
b}$ is a minimal binomial generator of $I_\L$. Since $|C_{\bf
b}|>2$ we obtain a contradiction.
\end{proof}

The following is now immediate:

\begin{cor1} If the generalized algebraic Scarf complex $({\bf
G_\L},\th)$ is a free resolution of $R/I_\L$ then $({\bf
G_\L},\th)$ is the algebraic Scarf complex and all $C$ that are
basic components of fibers are basic fibers themselves.
\end{cor1}

\section{Examples}
\label{examples_section}

In this section we compute the generalized algebraic Scarf complex
and compare it with the algebraic Scarf complex and the
generalized strongly algebraic Scarf complex for some examples. We
start with an example where the three complexes coincide.

\begin{ex1}\label{Scarf_resolution_not_generic} {\rm Let $\A$ be the semigroup
generated by the elements of the set $\{(4,0),(3,1), (1,3),
(0,4)\}$, and let $R=\Bbbk[a,b,c,d]$. The ideal $I_\A$ is
minimally generated by $bc-ad$, $ac^2-b^2d$, $b^3-a^2c$,
$c^3-bd^2$. The Scarf complex is a minimal resolution of
$R/I_{\A}$. We note that $I_\A$ is not a generic ideal in any
subring of $R$.}
\end{ex1}

In Example \ref{Koszul_example} we showed that when $I=\langle
f_1,\ldots, f_s \rangle$ is a complete intersection lattice ideal
where $f_i: i=1,\ldots, s$ is an $R$-sequence of binomials then
the Koszul complex on the $f_i$, $(\bf K_\bullet,\bf \theta)$,  is
a simple minimal resolution of $R/I_\L$.  In the following example
we compare $(\bf K_\bullet,\bf \theta)$ with the Scarf complex. We
also discuss  the largest index for which there is an
indispensable syzygy. First we give the following definition:

\begin{def1}
 We define {\em ideg($I_{\L}$)}, the indispensability degree of $I_{\L}$, to be the
 largest $t$ for
which there is an indispensable complex for $R/I_\L$ of length
$t$.
\end{def1}

\begin{ex1}{\rm Let the $f_i$ be as the above preceding remarks.
If the $f_i$ are indispensable generators for $I=\langle
f_1,\ldots, f_s \rangle$  then it is not hard to show that
$K_\bullet$ is a maximum indispensable complex $R/I_L$. In this
case ideg($I_{\L}$) is equal to the projective dimension of
$R/I_\L$. However the generalized algebraic Scarf complex may have
length just one, as is the case for the toric ideal $I=\langle
ae-fg, bd-cg\rangle$. If the $f_i$ are not indispensable then
ideg($I_{\L}$) can be 0, for example in the case of the toric
ideal $\langle x_1-x_2, x_2-x_3 \rangle$. }
\end{ex1}

The generalized algebraic Scarf complex is also computed in the
following two examples.

\begin{ex1} \label{strongly_algebraic_not
Scarf} {\rm Let $\A$ be the semigroup generated by the elements of
the set $\{(6,0),(4,2), (2,4), (0,6), (5,4)\}$, and let
$R=\Bbbk[a,b,c,d,e]$. Then $I_\A$ is minimally generated by
$-bc+ad, -b^2+ac, -c^2+bd, abd-e^2$. The corresponding lattice
$\L$ is given by the rows of the following matrix:
\[
 \begin{pmatrix} 1&-2&1&0&0\cr
                0&1&-2&1&0\cr
                0&2&1&0&-2\end{pmatrix}
                \]
A minimal resolution of $R/I_\A$, by \cite{CoCoA}, is of the form:
\[ 0\rightarrow R^2\rightarrow R^5\rightarrow R^4\rightarrow
R\rightarrow R/I_{\A}\rightarrow 0\] and the $i$-Betti degrees are
as follows:
\begin{itemize}
\item{} for $i=1$: $(6,6)$, $(8,4)$, $(4,8)$, $(10,8)$ \item{} for
$i=2$: $(8,10)$, $(10,8)$, $(14,16)$,
    $(16,14)$, $(18,12)$
\item{} for $i=3$: $(18,18)$, $(20,16)$.
\end{itemize}

For ${\bf b}= (6,6)$, the fiber $C_{\bf b}$ contains exactly 2
monomials: $bc, ad$. Their greatest common divisor is 1 and thus
the corresponding binomial is an indispensable generator of
$I_\A$. This is also the case for  ${\bf b}=(8,4)$,  and for ${\bf
b}=(4,8)$. For ${\bf b}=(10,8)$ the fiber $C_{\bf b}$ consists of
4 monomials, $ac^2, abd, b^2c, e^2$. Thus $\D_{\gcd}({\bf b})$
consists of a triangle and a point:

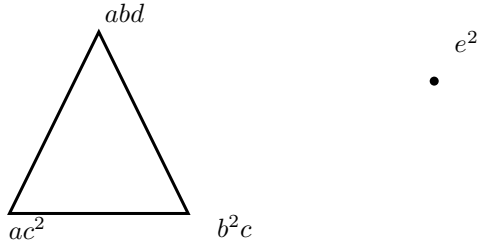
\begin{center}
\begin{pspicture}(0,-1.6129688)(6.7828126,1.6129688)
\pstriangle[linewidth=0.04,dimen=outer](1.4209375,-1.2654687)(2.44,2.48)
\psdots[dotsize=0.12](5.8809376,0.51453125)
\usefont{T1}{ptm}{m}{n}
\rput(1.7623438,1.4245312){$abd$}
\usefont{T1}{ptm}{m}{n}
\rput(0.47234374,-1.4354688){$ac^2$}
\usefont{T1}{ptm}{m}{n}
\rput(3.2223437,-1.4154687){$b^2c$}
\usefont{T1}{ptm}{m}{n}
\rput(6.3123436,1.0445312){$e^2$}
\end{pspicture}
\end{center}

\noindent The basic component of $\D_{\gcd}({\bf b})$, the
triangle, equals
 $C_{J_1}$
where $J_1=\{ (0,0,0,0,0)$, $(0,1,-2,1,0)$, $(1,-1,-1,0,0)\}$. We
note that $\bmax(J_1)$ $=\bmax(J)$ where $J=J_1\cup
\{(1,1,0,1,-2)\}$ and that $C_{\bf b}=C_{J}$. It follows that
$I_\A$ has a generator of $\A$-degree ${\bf b}$: any binomial
formed by taking the difference of $e^2$ from any element of
$C_{J_1}$ will do. Clearly this generator is not indispensable. It
is also immediate that $C_{J_1}$ is a basic component of $C_{\bf
b}$. Since $(10,8)$ is a minimal $1$-Betti degree it follows that
the image of $E_{C_{J_1}}$ gives a strongly indispensable syzygy.
Finally the fiber for $(8,10)$ is a basic fiber of cardinality 3
and equals $C_{J_2}$ for $J_2=\{ (0,0,0,0,0), (1,-1,-1,1,0),
(1,-2,1,0,0)\}$. The generalized algebraic Scarf complex equals
\[ 0\rightarrow R E_{C_{J_1}} \oplus R E_{C_{J_2}} \rightarrow R ^3
\rightarrow R\] and is strongly indispensable. We note that the
generalized algebraic Scarf complex equals the generalized
strongly algebraic Scarf complex and differs from the algebraic
Scarf complex. }
\end{ex1}

In the final example we give an ideal $I_\L$ for which the
generalized algebraic Scarf complex is not strongly indispensable.
In fact we give a fiber which  consists of two basic components.

\begin{ex1} \label{strongly_algebraic_not
Scarf} {\rm Let $\A$ be the semigroup of $\mathbb{Z}$ generated by
the 6 elements $3\cdot 13$, $4\cdot 13$, $5\cdot 13$, $3\cdot 14$,
$4\cdot 14$, $5\cdot 14$. In the ring $R=\Bbbk[a,\ldots, f]$, the
ideal $I_\A$ is generated by the binomials $-b^2+ac$, $e^2-df$,
$-a^3+bc$, $d^3-ef$, $-a^2b+c^2$, $d^2e-f^2$, and $bc^2-f^2d$. The
first 6 of these generators are indispensable binomials. The last
generator has $\A$-degree ${\bf b}=13\cdot 14$. In this case
$\D_{\gcd}({\bf b})$ has two basic components each of cardinality
3:

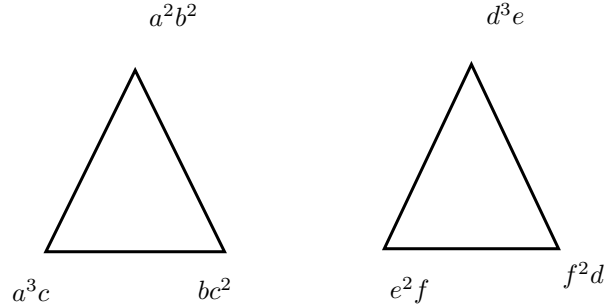
\begin{center}{
\begin{pspicture}(0,-2.0229688)(8.262813,2.0229688)
\pstriangle[linewidth=0.04,dimen=outer](1.8609375,-1.3554688)(2.44,2.48)
\pstriangle[linewidth=0.04,dimen=outer](6.3309374,-1.3154688)(2.38,2.52)
\usefont{T1}{ptm}{m}{n}
\rput(2.3723438,1.8145312){$a^2b^2$}
\usefont{T1}{ptm}{m}{n}
\rput(0.47234374,-1.8454688){$a^3c$}
\usefont{T1}{ptm}{m}{n}
\rput(2.9223437,-1.8254688){$bc^2$}
\usefont{T1}{ptm}{m}{n}
\rput(5.512344,-1.8454688){$e^2f$}
\usefont{T1}{ptm}{m}{n}
\rput(7.815,-1.6454687){$f^2d$}
\usefont{T1}{ptm}{m}{n}
\rput(6.782344,1.8345313){$d^3e$}
\end{pspicture}}
\end{center}

Thus any binomial formed by choosing one monomial from each
component is part of some minimal binomial generating set  of
$I_\A$. At the same time  there are two indispensable $2$-syzygies
of $\A$-degree $13\cdot 14$. Therefore these syzygies are not
strongly indispensable. Moreover we point out that there are
exactly two basic fibers with 3 elements: $C_{\bf 13\cdot 13}$ and
$C_{\bf 14\cdot 14}$. This means that no fiber has a basic
component of cardinality 4, since any such basic component would
imply the existence of 4 smaller basic fibers. A minimal
resolution of $R/I_\A$, by \cite{CoCoA}, is of the form:
\[ 0\rightarrow R^{4}\rightarrow R^{16}
\rightarrow R^{25}\rightarrow R^{19}\rightarrow R^7\rightarrow
R\rightarrow R/I_\A\rightarrow 0\ ,  \]  while the generalized
algebraic Scarf complex is of the form
\[ 0 \rightarrow R^{4}\rightarrow R^{6}
\rightarrow R\ .  \] We remark that the generalized strongly
algebraic Scarf complex is equal to the  algebraic Scarf complex
and differs from the generalized algebraic Scarf complex.}
\end{ex1}

As a final note we remark that it is easy to produce examples
where the three complexes are different, by combining the two
previous examples and using the technique of gluing semigroups,
see \cite{R}. \\

We finish with a list of related open questions.

\begin{itemize}

\item{} Determine all lattice ideals $I_\L$ so that the algebraic
Scarf complex is a minimal resolution of $R/I_\L$.

\item{} Determine  the maximum (strongly) indispensable  complex
for $R/I_\L$.

\item{} Determine all lattice ideals $I_\L$  so that the maximum
(strongly) indispensable complex is a minimal resolution of
$R/I_\L$.

\item{} Does there exist a $t$  such that if the monomial support
of all  i-syzygies of $R/I_\L$ for $i\le t$ form basic fibers,
then the Scarf complex is a free resolution for $I_\L$?

\item{} Does there exist a $t$ such that if ideg($I_{\L})\ge t$
then ideg($I_{\L}$) is equal to the projective dimension of
$R/I_\L$ and  the maximum indispensable complex is a free
resolution of $R/I_\L$?

\end{itemize}
\bigskip
{\bf Acknowledgment}
\bigskip
\par The authors would like to thank Ezra Miller for  useful
comments on this manuscript.

\end{document}